\def\draft{n}
\documentclass[headings]{amsart}
\usepackage{fullpage,amssymb,epic,eepic,epsfig,amsbsy,amsmath}

\theoremstyle{plain}

\newtheorem{theorem}{Theorem}
\newtheorem{metatheorem}{Meta-Theorem}

\newtheorem{proposition}{Proposition}[section]
\newtheorem{lemma}[proposition]{Lemma}
\newtheorem{corollary}[proposition]{Corollary}

\theoremstyle{definition}
\newtheorem{definition}[proposition]{Definition}

\newtheorem{question}{Question}

\theoremstyle{remark}

\newtheorem{remark}[proposition]{Remark}

\def\printname#1{
        \if\draft y
                \smash{\makebox[0pt]{\hspace{-0.5in}
                        \raisebox{8pt}{\tt\tiny #1}}}
        \fi
}

\newcommand{\psdraw}[2]
         {\begin{array}{c} \hspace{-1.3mm}
        \raisebox{-4pt}{\epsfig{figure=draws/#1.eps,width=#2}}
        \hspace{-1.9mm}\end{array}}

\newlength{\standardunitlength}
\catcode`\@=11
\long\def\@makecaption#1#2{%
     \vskip 10pt

\setbox\@tempboxa\hbox{
       \small\sf{\bfcaptionfont #1. }\ignorespaces #2}%
     \ifdim \wd\@tempboxa >\captionwidth {%
         \rightskip=\@captionmargin\leftskip=\@captionmargin
         \unhbox\@tempboxa\par}%
       \else
         \hbox to\hsize{\hfil\box\@tempboxa\hfil}%
     \fi}
\font\bfcaptionfont=cmssbx10 scaled \magstephalf
\newdimen\@captionmargin\@captionmargin=2\parindent
\newdimen\captionwidth\captionwidth=\hsize
\catcode`\@=12

\def\lbl#1{\label{#1}\printname{#1}}

\def\BN{\mathbb N}
\def\BZ{\mathbb Z}

\def\BQ{\mathbb Q}

\def\BC{\mathbb C}

\def\A{\mathcal A}

\def\calS{\mathcal S}

\def\a{\alpha}

\def\Ga{\Gamma}

\def\la{\langle}
\def\ra{\rangle}

\newcommand\strutn[2]{{\,^{#1}\!\!\frown^{#2}}}

\def\e{\epsilon}
\def\Ga{\Gamma}

\def\b{\beta}

\def\sub{\subset}

\def\longto{\longrightarrow}

\def\Lhat{\hat{\Lambda}}
\def\Lhathol{\hat{\Lambda}^{\mathrm{hol}}}

\def\fg{\mathfrak{g}}

\def\calD{\mathcal{D}}
\def\calA{\mathcal{A}}
\def\span{\mathrm{span}}
\def\calL{\mathcal{L}}

\def\calK{\mathcal{K}}
\newcommand{\Grad}{\operatorname{Grad}}
\newcommand{\tZ}{\tilde Z}
\def\ihs{integral homology sphere}
\def\tl{\mathrm{tl}}

\newwrite\transout
\openout\transout=figs.list

\newcommand{\ve}{\varepsilon}

\begin{document}


\title[Gevrey series in quantum topology]{
Gevrey series in quantum topology}

\author{Stavros Garoufalidis}
\address{School of Mathematics \\
         Georgia Institute of Technology \\
         Atlanta, GA 30332-0160, USA}
\email{stavros@math.gatech.edu,
URL: {\tt http://www.math.gatech.edu/$\sim$stavros }}

\author{Thang T.Q. Le}
\address{School of Mathematics \\
         Georgia Institute of Technology \\
         Atlanta, GA 30332-0160, USA}
\email{letu@math.gatech.edu,
URL: {\tt http://www.math.gatech.edu/$\sim$letu }}

\thanks{The authors were supported in part by the National
Science Foundation. \\
\newline
1991 {\em Mathematics Classification.} Primary 57N10. Secondary 57M25.
}

\date{February 12, 2007}


\begin{abstract}
Our aim is to prove that two formal power series of importance to quantum
topology are Gevrey.
These series are the Kashaev invariant of a knot (reformulated by Huynh and
the second author) and the Gromov norm of the LMO
of an integral homology 3-sphere. It follows that the power series associated
to a simple Lie algebra and a homology sphere is Gevrey.
Contrary to the case of analysis,
our formal power series are not solutions to differential equations with
polynomial coefficients.
The first author has conjectured (and in some cases proved, in joint work with
Costin) that our formal power series have resurgent Borel transform,
with geometrically interesting set of singularities.
\end{abstract}

\maketitle

\tableofcontents



\section{Introduction}
\lbl{sec.intro}

\subsection{Gevrey series}
\lbl{sub.gevrey}

A formal power series
\begin{equation*}
\lbl{eq.fx}
f(x)=\sum_{n=0}^\infty a_n \frac{1}{x^n} \in \BC[[1/x]]
\end{equation*}
is called {\em Gevrey-s} if there exists a positive constant $C$, such that
$$
|a_n| \leq  C^n n!^s
$$
for all $n>0$. Here, $x$ is supposed to be large. In other words, we will
order power series so that $1/x^n \gg 1/x^m$ iff $0 \leq n<m $.
Gevrey-0 series are well known: they are precisely the
convergent power series for $x$ in a neighborhood of infinity.
We will abbreviate Gevrey-1 by Gevrey. For example,
the following series
\begin{equation}
\lbl{eq.n!}
\sum_{n=0}^\infty n! \frac{1}{x^{n+1}}
\end{equation}
is Gevrey. Typically, Gevrey power series are divergent (for $x$ in a
neighborhood of infinity), and developing a meaningful calculus of
Gevrey power series is a well-studied subject; see \cite{Ha, Ra, Ec} and also
\cite{Ba}.
Gevrey power series appear naturally as formal power series solutions
to differential equations--linear or not. For example, the unique
formal power series solution to {\em Euler's equation}:
\begin{equation}
\lbl{eq.ode}
f'(x)+ f(x)=\frac{1}{x}
\end{equation}
is the series of Equation \eqref{eq.n!}. One can construct actual solutions of
the ODE \eqref{eq.ode} by suitably resumming the factorially divergent
series \eqref{eq.n!}, resulting in analytic functions with an essential
singularity at infinity; see \cite{Ha,Ra,Ec}.
The {\em resummation process}
of a Gevrey formal power series $f(x) \in \BQ[[1/x]]$
as in  \eqref{eq.fx} consists of the following steps:
\begin{itemize}
\item
consider its Borel transform $G(p)$, defined by:
\begin{equation*}
\lbl{eq.borelt}
G(p)=\sum_{n=1}^\infty a_n \frac{p^{n-1}}{(n-1)!} \in \BC[[p]]
\end{equation*}
Since $f(x)$ is Gevrey, it follows that $G(p)$ is analytic in a neighborhood
of $p=0$
\item
endless analytically continue $G(p)$ to a so-called resurgent
function,
\item
medianize if needed, and
\item
define the Laplace transform of $G(p)$ by:
\begin{equation*}
\lbl{eq.laplacet}
(\calL G)(x)=\int_0^\infty e^{-xp} G(p) dp
\end{equation*}
\end{itemize}
In the example the power series of \eqref{eq.n!}, its Borel transform $G(p)$
is given by:
$$
G(p)=\frac{1}{1-p}
$$
which is a resurgent (in fact, meromorphic) function with a single singularity
at $p=1$.

In general,
the output of a resummation is an analytic function (defined at least in
a right half-plane), constructed in a canonical way from the divergent
formal power series $f(x)$. In analysis, the resummation process commutes
with differentiation, and as a result one constructs actual solutions
of differential equations which are asymptotic to the formal power series
that one starts with.

A side corollary of resurgence (of importance to quantum
topology) is the existence of an asymptotic expansion of the coefficients
of the power series $f(x)$.  For a thorough discussion and examples, see
\cite{CG1}.

The above description highlights the necessity of
the Gevrey property, as a starting point of the
resummation.

\subsection{Formal power series in quantum topology}
\lbl{sub.qtop}

As mentioned before, a usual source of Gevrey series is a  differential
equation or a fixed-point problem.
Quantum topology offers a different source of Gevrey series
that do not seem to come from differential equations with polynomial
coefficients, due to the different structure of singularities of their
Borel transforms.
For example (and getting a little ahead of us),
the Kashaev invariant of two simplest knots (the trefoil ($3_1$),
and the figure eight ($4_1$)) are the power series:

\begin{eqnarray}
\lbl{eq.31}
F_{3_1}(x) &= & \sum_{n=0}^\infty (e^{1/x})_n \\
\lbl{eq.41}
F_{4_1}(x) &= & \sum_{n=0}^\infty (e^{1/x})_n (e^{-1/x})_n
\end{eqnarray}
where

\begin{equation*}
\lbl{eq.qq}
(q)_n=(1-q)\dots(1-q^n)
\end{equation*}
Notice that $(e^{1/x})_n \in 1/x^n \BQ[[1/x]]$, thus the power series
$F_{3_1}(x)$ and $F_{4_1}(x)$ are well-defined elements of
the formal power series ring $\BQ[[1/x]]$.

The power series $F_{3_1}(x)$ is the {\em Kontsevich-Zagier} power
series that was studied extensively by Zagier in \cite{Za}, and was
identified with the Kashaev invariant of the trefoil by Huynh and
the second author in \cite{HL}. In \cite{CG1}, Costin and the first
author gave an explicit formula for the Borel transform of
$F_{3_1}(x)$:

\begin{theorem}
\lbl{thm.CG1}\cite{CG1}
If $H_{3_1}(p)$ denotes the Borel transform of $e^{-1/(24x)} F_{3_1}(x)$, then
we have:
\begin{equation*}
\lbl{eq.Gformula}
H_{3_1}(p)=54 \sqrt{3} \pi
\sum_{n=1}^\infty \frac{\chi(n) n}{(-6 p + n^2 \pi^2)^{5/2}}.
\end{equation*}
where
\begin{equation*}
\lbl{eq.chi}
\chi(n)=
\begin{cases}
1 & \text{if} \,\, n \equiv 1,11 \bmod 12 \\
-1 & \text{if} \,\, n \equiv 5,7 \bmod 12 \\
0 & \text{otherwise.}
\end{cases}
\end{equation*}
\end{theorem}

Among other things, the above formula implies resurgence of the Borel
transform of the series $F_{3_1}(x)$ and locates explicitly the position
and shape of its singularities.

In \cite{CG2}, Costin and the first author prove by an abstract argument that
the Borel transform of the power series $F_{3_1}(x)$ and $F_{4_1}(x)$ are
resurgent functions.

The paper is concerned with two formal power series of importance to
quantum topology:
\begin{itemize}
\item
the Kashaev invariant of a knot,
\item
the LMO invariant of a closed 3-manifold.
\end{itemize}

Our aim is to prove that these series are Gevrey.

\subsection{The Gromov norm of the LMO invariant is Gevrey}
\lbl{sub.LMO}

Let us give a first impression the LMO invariant of
Le-Murakami-Ohtsuki, \cite{LMO}. It takes values in a (completed
graded) vector space $\A(\emptyset)$ of trivalent graphs, modulo
some linear AS and IHX relations:

\begin{equation*}
\lbl{eq.LMO} Z: \text{3-manifolds} \longto \calA(\emptyset)
\end{equation*}

The LMO invariant gives a meaningful definition to Chern-Simons
perturbation theory near a trivial flat connection. This is explained
in detail in \cite[Part I]{BGRT}. The trivalent
graphs are the Feynman diagrams of a $\phi^3$-theory (such as the
Chern-Simons theory) and their AS and IHX relations are
diagrammatic versions of the antisymmetry and the Jacobi identity of
the Lie bracket of a metrized Lie algebra.

The vector space $\A(\emptyset)$ has a grading (or degree) defined
by half the number of vertices of the trivalent graphs. Let
$\A_n(\emptyset)$ denote the subspace of $\calA(\emptyset)$ of
degree $n$.

As we discussed above, the LMO invariant takes values in $\calA(\emptyset)$.
In order to make sense of its Gevrey property, we need to replace
$\calA(\emptyset)$ by $\BQ[[1/x]]$. This is exactly what {\em weight systems}
do: they convert trivalent graphs into numerical constants; see \cite{B-N1}.
More precisely,
given a {\em simple Lie algebra} $\fg$, one can define a weight system map
(see \cite{B-N1}):

\begin{equation*}
W_{\fg}: \calA(\emptyset) \longto \BQ[[1/x]],
\end{equation*}
where each graph of degree $n$ is mapped into a rational number
times $1/x^n$. Combining the LMO invariant of a closed 3-manifold
$M$ with the weight system of a simple Lie algebra $\fg$, one gets a
formal power series:
\begin{equation}
\lbl{eq.FSdefg}
F_{\fg,M}(x)= W_{\fg}(Z_M) \in \BQ[[1/x]]
\end{equation}
This power series is equal to the  Ohtsuki series, defined by
Ohtsuki \cite{Oh1} for $\fg = sl_2$ and then by the second author for
all simple Lie algebras \cite{Le5}. As nice as weight systems are,
the power series still depends on Lie algebras; moreover it is known
that not all weight systems come from Lie algebras, \cite{Vo}.

Ideally, we would like to replace the graph-valued invariant $Z_M
\in \calA(\emptyset)$ by a single series $|Z_M| \in \BQ[[1/x]]$ so
that

\begin{itemize}
\item[(a)]
$|Z_M| \in \BQ[[1/x]]$ is Gevrey,
\item[(b)]
The Gevrey property of $|Z_M|$ implies the Gevrey property of
$F_{\fg,M}(x)$ for all simple Lie algebras,
\item[(c)]
$|Z_M| \neq 1$ iff $Z_M \neq 1$.
\end{itemize}

Can we accomplish this at once? A simple idea, the Gromov norm,
allows us to achieve this.

\begin{definition}
\lbl{def.gnorm1}
Consider a  vector space $V$ with a subset $b$ that
spans $V$.  For $v\in V$, define $b$-norm by
$$
|v|_{b}=\inf{\sum_j |c_j|}
$$
where the infimum is taken over all presentations of the form $v
=\sum_j c_jv_j, v_j \in b$.
\end{definition}

For example, consider $V = \BQ[q^{\pm 1}]$ -- the space of Laurent
polynomials in $q$ with rational coefficients,
and $b$ the set $\{q^n\,|\, n \in \BZ\}$. In this
case the norm of a Laurent polynomial $f(q)$ is known as its
$l^1$-norm, denoted by $||f(q)||_1$.

We will apply Definition \eqref{def.gnorm1} to $V= \A(\emptyset)$,
with $b$ is the set of {\em trivalent graphs},
and we will denote $|v|_b$ simply by $|v|$ for $v \in \A(\emptyset)$.
For a precise definition of what is a trivalent graph, see Section
\ref{sec.LMOgevrey}.

\begin{definition}
\lbl{def.gnorm}
\rm{(a)}
For an element $v \in \A(\emptyset)$ let
$\Grad_n(v)$ be the part of degree $n$ of $v$. The Gromov norm of
$v$ is defined as
\begin{equation*}
\lbl{eq.gnorm} |v|=\sum_{n=0}^\infty |\Grad_n (v)| \frac{1}{x^n} \in
\BQ[[1/x]].
\end{equation*}
\rm{(b)} Let us say that $v \in \A(\emptyset)$ is {\em Gevrey-$s$} iff
$|v| \in \BQ[[1/x]]$ is Gevrey-$s$.
\end{definition}
It is easy to see that $|v|=1$ iff $v=1$.
Here, $1 \in \A(\emptyset)$ denotes the element of degree $0$
which is $1$ times the empty trivalent graph.

Our next theorem explains a Gevrey property of the LMO invariant.

\begin{theorem}
\lbl{thm.gLMO} For every \ihs\ $M$, $|Z_M| \in
\BQ[[1/x]]$ is Gevrey.
\newline
Moreover, $|Z_M|=1$ iff $Z_M=1$.
\end{theorem}
Theorem \ref{thm.gLMO} and an easy estimate implies the following:

\begin{theorem}
\lbl{thm.4} For every closed 3-manifold and every simple Lie algebra
$\fg$, the Ohtsuki series $F_{\fg,M}(x)$ is Gevrey.
\end{theorem}

A key ingredient in the definition of the LMO invariant is the {\em
Kontsevich integral} $Z_L$ of a framed link $L$ in $S^3$. The
Gromov norm of $Z_L$ can be defined in a similar fashion, see
Section \ref{sec.LMOgevrey}. In fact, Theorem \ref{thm.gLMO}
motivates (and even requires) to consider the Gromov norm $|Z_L|
\in \BQ[[1/x]]$ of the Kontsevich integral.

\begin{theorem}
\lbl{thm.gknot} For every framed link $L$ in $S^3$, $|Z_L| \in
\BQ[[1/x]]$ is Gevrey-0.
\end{theorem}
Recall that a power series $f(x)$ is Gevrey-0 iff $f(x)$ is a convergent
power series for $x$ near $\infty$.

Theorems \ref{thm.gLMO} and \ref{thm.gknot}
are a special case of the following guiding principle, which
we state as a Meta-Theorem:

\begin{metatheorem}
Asymptotic power series that appear in constructive quantum
field theory (and in particular, in Quantum Topology) are resurgent
functions--and in particular, Gevrey of some order
(usually, order $1$).
\end{metatheorem}

Let us comment that
the factorial growth of power series in perturbative quantum field
theory is usually due to the factorial growth of the number of Feynman
diagrams; see for example Lemma \ref{lem.normalized}.
The contribution of each Feynman diagram is growing exponentially only;
see for example Lemma \ref{lem.weightbound}.

\subsection{The Habiro ring}
\lbl{sub.ring}

So far we discussed how perturbative quantum field theory leads to
Gevrey power series \eqref{eq.FSdefg}. Examples of such series
(for knots, rather than 3-manifolds) were given in Equations
\eqref{eq.31} and \eqref{eq.41}.

In the remaining of this section, we will concentrate with the case
of $\fg=\mathfrak{sl}_2$. Our aim is to give a non-perturbative
explanation of the Gevrey property of the power series
$F_{\mathfrak{sl}_2,M}(x)$, which we will abbreviate by $F_M(x)$ in
this section. In fact, we will be dealing with a formal power series
invariant of knotted objects:

\begin{equation}
\lbl{eq.FS}
F: \text{Knotted Objects} \longto \BQ[[1/x]]
\end{equation}
where a {\em knotted object} (denoted in general by $\calK$) will be
either a knot $K$ in 3-space or an \ihs\ $M$. We
already discussed the series $F_M(x):=F_{\mathfrak{sl}_2,M}(x)$. In
the case of a knot $K$, the power series $F_K(x)$ will be defined
below.

In the absence of a rule (such as a differential equation) for the
power series $F_{\calK}(x)$, or an explicit formula (in the style of
\eqref{eq.31} or \eqref{eq.41}), how can one prove that our power
series are Gevrey? It turns out that the power series $F_{\calK}(x)$
have a certain ``shape'' which explains their Gevrey (and
conjectural resurgence) property. Such a shape was discovered by
Habiro, who considered the {\em cyclotomic completion} of the ring
of Laurent polynomial (the so-called {\em Habiro ring})
\begin{equation*}
 \Lhat=\lim_{\leftarrow n} \BZ[q^{\pm 1}]/( (q)_n )
\end{equation*}
As a set, it follows that the Habiro ring is:
\begin{equation*}
\Lhat=\{ f(q)=\sum_{n=0}^\infty f_n(q) (q)_n \,| \, f_n(q) \in
\BZ[q^{\pm 1}] \}
\end{equation*}
Habiro showed a number of key properties of the ring $\Lhat$; see
\cite{H2}. For our purposes, it will be important that elements
$f(q)$ of the Habiro ring have Taylor series expansions at $q=1$,
and that they are uniquely determined by their Taylor series. In
other words, the map from $\Lhat$ to $\BZ[[q-1]]$, sending $f(q)$ to
its Taylor series at 1, is injective. Let

\begin{equation*}
 T: \Lhat \longto   \BQ[[1/x]]
\end{equation*}

be the map defined so that $T(f(q))$ is obtained from the Taylor
series of $f(q)$ by the substitution $q= e^{1/x}$. Then $T$ is
injective.

 In the case of an \ihs\ $M$, Habiro
proved that the series $F_M(x)$ comes from a (unique) element
$\Phi_M(q)$ of the Habiro ring. In the case of a knot, Huynh and the
second author observe in \cite{HL} that the {\em Kashaev invariant}
of a knot $K$ also comes from an element $\Phi_K(q)$ of the Habiro
ring. In that case, we define $F_K(x)=(T
\Phi_K)(x)=\Phi_K(e^{1/x})$.

In other words, we have a map

\begin{equation*}
\Phi: \text{Knotted Objects} \longto \Lhat
\end{equation*}
such that
\begin{equation*}
F=T \circ \Phi.
\end{equation*}
Thus, instead of writing
$$
F_{\calK}(x)=\sum_{n=0}^\infty a_{\calK,n} \frac{1}{x^n}
$$
for $a_{\calK,n} \in \BQ$, we may write:
$$
F_{\calK}(x)=\sum_{n=0}^\infty f_{\calK,n}(e^{1/x}) (e^{1/x})_n
$$
for suitable polynomials $f_{\calK,n}(q) \in \BZ[q^{\pm 1}]$. Keep in mind that
the polynomials $f_{\calK,n}(q)$ are not unique. For example, we have
the following identity in the Habiro ring:
\begin{equation*}
\sum_{n=0}^\infty -q^{n+1} (q)_n
=1.
\end{equation*}
Most importantly for us, without
any additional information about the polynomials $f_{\calK,n}(q)$ one cannot
expect that the series $F_{\calK}(x)$ is Gevrey.
The information can be formalized by introducing two subrings
of $\Lhat$. We need an auxiliary definition.

\begin{definition}
\lbl{def.qholo}
\rm{(a)}
We say that a sequence $(f_n(q))$ of Laurent polynomials is
{\em $q$-holonomic}
if it satisfies a linear $q$-difference equation of the form:
$$
a_d(q^n,q) f_{n+d}(q) + \dots  a_0(q^n,q) f_{n}(q)=0
$$
for all $n \in \BN$, where $a_j(u,v) \in \BZ[u^{\pm 1},v^{\pm 1}]$
for $j=0,\dots,d$ and $a_d \neq 0$.
\newline
\rm{(b)} We say that a sequence $(f_n(q))$ of Laurent polynomials is
{\em nicely bounded} if there exist the bounds on their span and
coefficients: There are constants $C, C'>0$ such that for $n>0$,
\begin{eqnarray}
\lbl{eq.degbound}
\span_q f_n(q) & \sub & [-C' \, n^2, C' \, n^2] \\
\lbl{eq.coeffbound}
||f_n(q)||_1 & \leq & C^n.
\end{eqnarray}
\end{definition}

Now, we may define the following subrings of the Habiro ring.

\begin{definition}
\lbl{def.qholob}
\rm{(a)}
We define:
\begin{equation*}
\Lhathol=\{
f(q)=\sum_{n=0}^\infty f_n(q) (q)_n
\,| f_n(q) \in \BZ[q^{\pm 1}],
\,\, (f_n(q)) \,\, \text{is $q$-holonomic}
\}
\end{equation*}
\newline
\rm{(b)}
We define:
\begin{equation*}
\Lhat^b=\{
f(q)=\sum_{n=0}^\infty f_n(q) (q)_n
\,| f_n(q) \in \BZ[q^{\pm 1}], \,\, (f_n(q)) \,\, \,\, \text{is nicely
bounded} \}
\end{equation*}
\end{definition}

It is easy to see that $\Lhat^h$ and $\Lhat^b$ are subrings of
$\Lhat$. Observe that $\Lhat^h$ is a countable ring, whereas $\Lhat$
and $\Lhat^b$ are not.

It is easy to show that if $(f_n(q))$ is a $q$-holonomic sequence of
Laurent polynomials, then it satisfies \eqref{eq.degbound}. On the
other hand, the authors do not know the answer to the following
question.

\begin{question}
\lbl{que.bh}
Is it true that $\Lhat^h$ is a subring of $\Lhat^b$?
\end{question}

\subsection{Gevrey series from the Habiro ring}
\lbl{sub.results}

Independently of the answer to the above question, we have:

\begin{theorem}
\lbl{thm.1}
For every knotted object $\calK$ we have:
\begin{equation}
\lbl{eq.IMhb}
\Phi_{\calK}(q) \in \Lhathol \cap \Lhat^b.
\end{equation}
\end{theorem}

Our next theorem relates the ring $\Lhat^b$ with Gevrey series.


\begin{theorem}
\lbl{thm.2}
If $f(q) \in 
\Lhat^b$, then $f(e^{1/x})\in \BQ[[1/x]]$ is Gevrey. 
\end{theorem}

Theorems \ref{thm.1} and \ref{thm.2} imply the promised result.

\begin{theorem}
\lbl{thm.3}
For every knotted object $\calK$, the power series $F_{\calK}(x)$
is Gevrey.
\end{theorem}

If $M$ is an \ihs\, then the above theorem gives an
independent proof that the series $F_M(x)$ is Gevrey.

\subsection{What next?}
\lbl{sub.whatnext}

As was mentioned in Section \ref{sub.gevrey}, a Gevrey series is the
input of a resummation process. In \cite{CG2} we conjecture that the
series $F_{\calK}(x)$ of every knotted object $\calK$ can be
resummed. In other words, we conjecture that the Borel  transform
$G_{\calK}(p)$ of $F_{\calK}(x)$ is a resurgent function, with
singularities given by geometric invariants of the knotted object
$\calK$. This conjecture is true for the two simplest knots $3_1$
and $4_1$ and for several elements of $\Lhathol$; see \cite{CG1} and
\cite{CG2}. Based on this partial evidence, we pose the following
questions:

\begin{question}
\lbl{que.1}
If $f(q) \in \Lhathol$, is it true that its Taylor series $(Tf)(x) \in
\BQ[[1/x]]$ has resurgent Borel transform?
\end{question}

\begin{question}
\lbl{que.2} Is it true that the Gromov norm $|Z_M| \in \BQ[[1/x]]$
of of the LMO invariant of an \ihs\ has resurgent
Borel transform?
\end{question}

\begin{question}
\lbl{que.3} Is it true that the Gromov norm $|Z_L| \in \BQ[[1/x]]$
of the Kontsevich integral of a framed link in $S^3$ is a resurgent
function?
\end{question}

\subsection{Plan of the proof}
\lbl{sub.plan}

Since Gevrey series are not familiar objects in quantum topology,
we have made an effort to motivate their appearance and usefulness
in quantum topology. For the analyst, we would like to point
out that our Gevrey series (and their expected resurgence properties)
are not expected to be solutions of differential equations (linear or not)
with polynomial coefficients. Thus, our results are new from this perspective.

We have also separated into different sections results from quantum
topology and from asymptotics.

In Section \ref{sec.LMOgevrey} we discuss in detail the LMO invariant,
starting from the necessary discussion of the Kontsevich integral
of a framed link in 3-space. Basically, the LMO invariant of a 3-manifold
is obtained by the (suitably normalized) Kontsevich integral of a surgery
presentation link, after we glue all legs.
We will use combinatorial counting arguments to bound the number
of unitrivalent graphs,  as well as the original definition of
the Kontsevich integral to estimate the coefficients
of these graphs, before and after the gluing of the legs.
In addition in Section \ref{sub.ohtsukiseries}
we show that various analytic reparametrizations of the
LMO invariant (such as the Ohtsuki series) are Gevrey.
This ends the perturbative quantum field theory discussion of the paper.

In Section \ref{sec.proofthm1} we give a nonperturbative explanation
of the Gevrey property of our power series for the simple Lie algebra
$\mathfrak{sl}_2$. In that case, the Kontsevich integral is replaced by
the colored Jones function of a link. The latter is a multisequence of
Laurent polynomials. We discuss two key properties of the
colored Jones function: $q$-holonomicity (introduced in \cite{GL1})
and integrality, introduced by Habiro in \cite{H1,H2}. Together with Habiro's
definition of $\Phi_M(q)$ (given in terms of a surgery presentation of $M$),
$q$-holonomicity implies that $\Phi_M(q) \in \Lhathol$, and integrality
implies that $\Phi_M(q) \in \Lhat^b$. Combined together with
Theorem \ref{thm.2} (shown in the next section), they give a proof of
Theorem \ref{thm.3}.

Finally,
in Section \ref{sec.main} we use elementary estimates to give a proof
of Theorem \ref{thm.2}.

\subsection{Acknowledgement}
An early version of this paper was presented by lectures of the first author
in Paris VII in the summer of 2005.
The first author wishes to thank G. Masbaum and P. Vogel for their hospitality.

\section{The LMO invariant is Gevrey}
\lbl{sec.LMOgevrey}

In this section we will give a proof of Theorem \ref{thm.gLMO}.

Omitting technical details, the {\em Aarhus version} of the LMO
invariant \cite[Part II \& III]{BGRT} is defined as follows. Suppose
an integral homology $M$ is obtained from $S^3$ by surgery on a
framed link $L$.

\begin{itemize}
\item
Consider a  presentation of $L$ as the closure of a framed string
link $T$.
\item
Consider the suitably normalized {\em Kontsevich integral} $\check
Z_T$ of the string link $T$. It takes values in a completed
$\BQ$-vector space of vertex-oriented unitrivalent graphs (Jacobi
diagrams) with legs colored by  the components of $L$.
\item
Separating out the strut part from $\check Z_T$ and closing  we get
the formal Gaussian integral $\int \check Z_T$, which takes values
in the algebra $\A(\emptyset)$ of Jacobi diagrams without legs.
\item
Finally, normalize $\int \check Z_T$ in a minor way to get the LMO
invariant $Z_M$.
\end{itemize}

The precise definition will be recalled later. To prove Theorem
\ref{thm.gLMO} we will need to have an estimate
\begin{itemize}
\item[(a)]
for the norm of the Kontsevich integral and
\item[(b)] for the norms of the maps appearing in the
definition of the LMO invariant.
\end{itemize}
To get the desired estimates it will be simpler to exclude Jacobi diagrams
with tree components. This is guaranteed when $L$ is a boundary link.
And it suffices since every \ihs\ can be obtained by surgery along a
unit-framed boundary link.

\subsection{Jacobi diagrams}
\lbl{sub.jacobi}

We quickly recall the basic definitions
and properties here, referring the details to \cite{B-N1,BLT}.

\lbl{Preliminaries}
 An {\em open Jacobi diagram}  is a
vertex-oriented uni-trivalent graph, i.e., a graph with univalent
and trivalent vertices together with a cyclic ordering of the edges
incident to the trivalent vertices.   A univalent vertex is called
{\em a leg}, and trivalent vertex is also called an {\em internal
vertex}.
The {\em degree} of an open Jacobi diagram is half the number of
vertices (trivalent and univalent). The i-degree is the number of
internal vertices, and the e-degree is the number of legs.

Suppose  $X$ is a compact oriented 1-manifold (possibly with
boundary) and $Y$ a finite set. A {\em Jacobi diagram based on
$X\cup Y$} is a graph $D$ together with a decomposition $D = X \cup
\Gamma$, where $\Gamma$ is an open Jacobi diagram with some legs
labeled by elements of $Y$, such that $D$ is the result of gluing
all the non-labeled legs of $\Gamma$ to distinct interior points of
$X$. Note that repetition of labels is allowed. The {\em degree} of
$D$, by definition, is the degree of $\Gamma$.

The space $\A^f(X, Y)$ is the vector space over $\BQ$ spanned by
Jacobi diagrams based on $X\cup Y$ modulo the usual AS,
IHX and STU relations (see \cite{B-N1}). The
completion of $\A^f(X, Y)$ with respect to degree is denoted by
$\A(X, Y)$.

Of special interest are the following flavors of Jacobi diagrams:

\begin{itemize}
\item[(a)]
$(X,Y)=(\circlearrowleft_m,\emptyset)$, where
$X=\circlearrowleft_m$ the union of $m$ numbered, oriented
circles. Then,  $\A(X,Y)=
\A(\circlearrowleft_m)$ is the space where the
Kontsevich integral of $m$-component framed link lies.
\item[(b)]
$(X,Y)=(\uparrow_m,\emptyset)$, where  $X= \uparrow_m$, the union of
$m$ numbered, oriented intervals. Then,
$\A(X,Y)=\A(\uparrow_m)$ is the space where the
Kontsevich integral of $m$-component framed braid (or a dotted Morse link)
lies.
\item[(c)]
$(X,Y)=(\emptyset, \{1,2,\dots,m\})$. Then, we denote
$\A(X,Y)$ by $\A(\star_m)$
\item[(d)]
$(X,Y)=(\circlearrowleft_r \cup \uparrow_s,\emptyset)$. Then
$A(X,Y)$ is the space where the Kontsevich integral of a tangle $T$ lies,
where $r$ is the number of interval components of $T$ and
$s$ is the number of circle components of $T$.
\end{itemize}
The spaces $\A(\circlearrowleft_m)$, $\A(\uparrow_m)$
are related with an obvious projection map:

\begin{equation}
\lbl{eq.projp}
p: \A(\uparrow_m) \to \A(\circlearrowleft_m)
\end{equation}
which identifies the two end points of each interval in $\uparrow_m$.

The spaces $\A(\uparrow_m)$ and $\A(\star_m)$ are also
related with a {\em symmetrization map}

\begin{equation}
\lbl{eq.symmetrization}
\chi:\A(\star_m) \rightarrow \A(\uparrow_m)
\end{equation}
which is a linear map  defined on a diagram $\Gamma$ by taking the average over
all possible ways of ordering the legs labeled by $j$, $1\le j \le
m$, and attach them to the $j$-th oriented interval. It is known
that $\chi$ is a vector space isomorphism \cite{B-N1}.

\begin{remark}
\lbl{rem.2mult}
Note that $\A(\uparrow_m)$ is an algebra, where the product of two Jacobi
diagrams is obtained by placing (or stacking) the first on top of the second.
$\A(\star_m)$ is also an algebra, where the product of two diagrams is their
disjoint union. However, the map $\chi$, which is a vector space
isomorphism, is not an algebra isomorphism. To get
algebra isomorphism one needs the wheeling map, see \cite{BLT}. To
avoid confusion we use $\#$ to denote the product in
$\A(\uparrow_m)$ and $\sqcup$ the product in $\A(\star_m)$.
\end{remark}

The diagonal map $\Delta^{(m)} : \A(\star_1) \to \A(\star_m)$ is a
linear map defined on a Jacobi diagram $\Gamma\in \A(\star_1)$ by
taking the sum of all possible Jacobi diagrams $\Gamma'\in
\A(\star_m)$ such that if we switch all the labels  in $\Gamma'$ to
$1$, then we from $\Gamma'$ we get $\Gamma$. It is clear that if
$\Gamma$ has $k$ legs, then there are $m^k$ such $\Gamma'$.

Suppose $X=Y=\emptyset$. The space $\A(\emptyset)$ is the space in
which lie the values of the LMO invariants of 3-manifolds
\cite{LMO}. With disjoint union as the product, $\A(\emptyset)$
becomes a commutative algebra, and all other $\A(X, Y)$ have a
natural $\A(\emptyset)$-module structure.

An open Jacobi diagram is {\em tree-less} if none of its connected
components is a tree. Let $\A^{\tl}(\star_m)$ be the subspace of $\A(\star_m)$
spanned by treeless Jacobi diagrams. The following is obvious but
useful later.

\begin{lemma} Suppose $x\in \A^{\tl}(\star_m)$ has i-degree $n$,
then the e-degree of $x$ is less than or equal to $n$. \lbl{101}
\end{lemma}

\subsection{Norm of Jacobi diagrams}

The set of Jacobi diagrams based on $X\cup Y$ clearly spans the
space $\A(X,Y)$. The norm of $v\in \A(X,Y)$ with respect to this
spanning subset is denoted simply by $|v|$. We will say that
$v \in \A(X,Y)$ is Gevrey-$s$ iff $|v| \in \BQ[[1/x]]$ is Gevrey-$s$.

It is clear that if the product $vu$ can be defined, then $|vu| \le
|v| |u|$. Since the product of two Gevrey-$s$ power series is Gevrey-$s$,
and the inverse of a Gevrey-$s$ series with nonzero constant term
is Gevrey-$s$, it follows that:

\begin{lemma}
\lbl{205}
\rm{(a)} If $v,u\in A(X,Y)$ are Gevrey-$s$ and the product $vu$
can be defined, then $vu$ is also Gevrey-$s$.
\newline
\rm{(b)}If $v\in \A(\emptyset)$ has non-zero constant term and is
Gevrey-$s$, then $1/v\in \A(\emptyset)$ is Gevrey-1. 
\end{lemma}

For an element $v\in \A^{\tl}(\star_m) \subset \A(\star_m)$, in addition to the
above norm, there is another one defined using the spanning set of
{\em treeless Jacobi diagrams}.

\begin{lemma}
\lbl{lem.equalnorms}
The above two norms are equal.
\end{lemma}

\begin{proof}
The lemma follows at once from the fact that the subspace spanned by
Jacobi diagrams other than treeless ones intersects $\A^{\tl}(\star_m)$ only
by the zero vector.
\end{proof}

Recall the symmetrization map $\chi$ from \eqref{eq.symmetrization}.
The next proposition estimates the norm of $\chi$ and $\chi^{-1}$.

\begin{proposition}
\lbl{204}
\rm{(a)} For every $v\in \A(\star_m)$, one has $|\chi(v)| \le |v|$. In other
words, the operator $\chi$ has norm less than or equal to 1.
\newline
\rm{(b)} Suppose $x\in \A(\uparrow_m)$ has e-degree $k \ge 1$, then
$|\chi^{-1}(v)| \le 2k |v|$.
\newline
\rm{(c)} For any $v\in \A(\star_1)$ of e-degree $k$, one has
$|\Delta^{(m)}(v)| \le m^k |v|$.
\end{proposition}

\begin{proof}
(a) and (c) follows immediately from the definition.
We give here the proof of (b).

We use induction. Suppose the statement holds true when $v$ has
e-degree $<k$. It is enough to prove for the case when $v = \Gamma
\cup \uparrow_m$, where $\Gamma$ is an open Jacobi diagram with $k$
legs. Using the STU relation, one can see that $u:= \chi(\Gamma)-v$
has e-degree $< k$. One has

$$
|u| = | \chi(\Gamma)-v| \le |\chi(\Gamma)| + |v| \le 2.
$$
By induction, $|\chi^{-1}(u)| < 2(k-1)$. Since $v = u +
\chi(\Gamma)$, we have $\chi^{-1}(v) = \chi^{-1}(u) + \Gamma$, and
hence
$$
| \chi^{-1}(v) | \le |\chi^{-1}(u)| + |\Gamma| \le 2(k-1) +1 < 2k.
$$
\end{proof}

\subsection{The unknot}
\lbl{sub.unknot}

Let $w_{2n}\in \A(\star_1)$ be the {\em wheel} with $2n$ legs. It is
the open Jacobi diagram consisting of a circle and $2n$ intervals
attached to it. For example,
$$
w_4=\psdraw{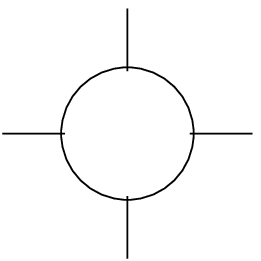}{0.3in}
$$
Define

\begin{equation*}
\nu=\exp \left( \sum_{n=1}^\infty b_{2n}\omega_{2n} \right) \in
\A(\star_1) \qquad \text{and} \qquad \sqrt \nu=\exp \left(
\frac{1}{2}\sum_{n=1}^\infty b_{2n}\omega_{2n} \right) \in
\A(\star_1).
\end{equation*}
where the {\em modified Bernoulli numbers} $b_{2n}$ are defined by the
power series expansion

\begin{equation}
\lbl{eq:MBNDefinition}
\sum_{n=0}^\infty b_{2n}x^{2n} = \frac{1}{2}\log\frac{\sinh x/2}{x/2}.
\end{equation}
Notice that
$$
\frac{1}{2}\log\frac{\sinh x/2}{x/2}=\frac{1}{48} x^2 -\frac{1}{5760} x^4
+\frac{1}{362880} x^6 + \dots
$$
The modified Bernoulli numbers are related to the
zeta function
$\zeta(n) := \sum_{k=1}^\infty k^{-n}$ by

\begin{equation}
\lbl{001}
b_{2n} = (-1)^{n-1} \frac{\zeta(2n)}{2n\, (2\pi)^{2n}}.
\end{equation}

In \cite{BLT} it was shown that $\chi^{-1}(\nu)$ is the Kontsevich
integral of the unknot.

\begin{proposition}
\lbl{002}
The series $\nu$ and $\sqrt \nu $ are Gevrey-0.
\end{proposition}

\begin{proof}
Since $|w_{2n}| \le 1$, it is enough to to show that the
series $\exp (\sum_n |b_{2n}| x^{-2n}) $ is convergent for large enough $x$.
Since $|b_{2n}|=(-1)^{n-1} b_{2n}$, it follows that
\begin{eqnarray*}
\exp(\sum_n |b_{2n}| x^{-2n}) &=&
\exp(-\sum_n b_{2n} (i/x)^{2n}) \\
&=&
\exp\left(-\frac{1}{2} \log \frac{\sin(1/(2x))}{1/(2x)} \right) \\
&=& \sqrt{\frac{1/(2x)}{\sin(1/(2x))}}
\end{eqnarray*}
The latter converges for $|x| > 1/(4 \pi)$.
\end{proof}

\begin{question}
\lbl{que.nu}
Is it true that
$$
|\nu|= \sqrt{\frac{1/(2x)}{\sin(1/(2x))}}?
$$
\end{question}

Since  $\Grad_n\nu$, the part of degree $n$  of $\nu$, has e-degree
$2n$, Proposition \ref{204}(c) implies that

\begin{corollary}
\lbl{207}
For every positive integer $m$ the series
$\Delta^{(m)}(\nu) \in \A(\star_m)$ is Gevrey-0.
\end{corollary}

\subsection{The Kontsevich integral}
\lbl{sub.knot}

The {\em framed} Kontsevich integral of a {\em framed} tangle $T$
takes value in $\A(T)$, see for example \cite{LM,B-N2,BLT}.
This is a slight modification of the original integral defined by Kontsevich
\cite{Ko}. The framed Kontsevich integral depends on the positions
of the boundary points. To get rid of this dependence one has to
choose standard positions for the boundary points. It turns out that
the best positions are in a limit, when all the boundary points
go to one fixed point. In addition one has to regularize the Kontsevich
integral in the limit. In the limit one has to keep track of the
order in which the boundary points go to the fixed point. This leads
to the notion of {\em parenthesized framed tangle}. The latter were
called  {\em $q$-tangle} in \cite{LM} and {\em non-commutative tangles}
in \cite{B-N2}. For details, see \cite{LM} and \cite{B-N2}.

In all framed tangles in this paper, we assume that a
non-associative structure is fixed. Theorem \ref{thm.gknot} is a
special case of the following theorem.

\begin{theorem}
\lbl{206}
For every framed tangle $T$ the Kontsevich integral
$Z_T$ is Gevrey-0.
\end{theorem}

\begin{proof}
This follows from the facts that
\begin{itemize}
\item[(a)] the Kontsevich integral satisfies a locality property.
In other words, a framed tangle is the assembly
(i.e., the product) of
elementary blocks of three kinds: local extrema, crossings, and
change of parenthetization; for a computerized example, see
\cite[Sec.1.2]{B-N2}.
The Kontsevich integral is the
corresponding product of the invariants of the elementary blocks,
\item[(b)] The invariants of each blocks are Gevrey-$0$, and
\item[(c)] The product of Gevrey-$0$ series is Gevrey-$0$ by Lemma
\ref{205}.
\end{itemize}
More precisely, the Kontsevich integral of the elementary blocks is given
by:
$$
\psdraw{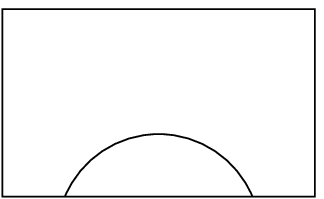}{0.7in}, \psdraw{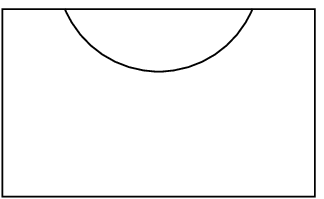}{0.5in} \to \sqrt{\nu}
$$
$$
\psdraw{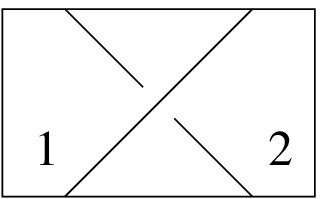}{0.7in} \to \exp\left(\frac{1}{2} \strutn{1}{2}\right)
$$
$$
\psdraw{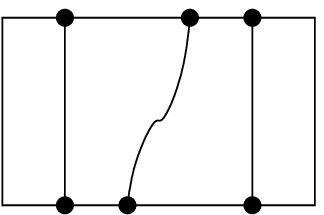}{0.7in} \to \Phi
$$
where the {\em strut} $\strutn{1}{2}$ is the only open Jacobi
diagram homeomorphic to an interval, and $\Phi$ is any {\em associator}.
For a definition of an associator, see \cite{Dr} and also \cite{B-N2,B-N3}.
Proposition \ref{002} implies that $\sqrt{\nu}$ is Gevrey-$0$.
In \cite{LM}, J. Murakami and the second author
gave an explicit formula for the {\em KZ-associator} $\Phi^{\rm{KZ}}$:
\begin{eqnarray}
\lbl{eq.PhiKZ}
\Phi^{\rm{KZ}} &=& 1+\sum_{l=1}^\infty
\sum_{\mathbf{a},\mathbf{b},\mathbf{p},\mathbf{q}}
(-1)^{|\mathbf{b}|+|\mathbf{p}|} \eta(\mathbf{a}+\mathbf{p},
\mathbf{b}+\mathbf{q})
\binom{\mathbf{a}+\mathbf{p} }{\mathbf{b}+\mathbf{q} }
B^{|\mathbf{q}|} (A,B)^{(\mathbf{a},\mathbf{b})} A^{|\mathbf{p}|} \\
\notag
&=& 1 + \frac{1}{24} [A,B] -\frac{\zeta(3)}{(2 \pi i)^3} [A,[A,B]]
+\dots
\end{eqnarray}
where
$$
A=\psdraw{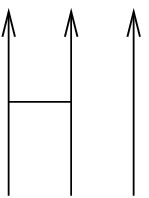}{0.3in}, \qquad B=\psdraw{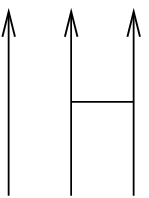}{0.3in},
$$
$[X,Y]=XY-XY$,
\begin{equation}
\zeta(a_1, a_2,\dots,a_k)=\sum_{n_1 < n_2 < \dots < n_k \in \BN}
n_1^{-a_1} n_2^{-a_2} \dots n_k^{-a_k}
\end{equation}
are the multiple zeta numbers and for $\mathbf{a}=(a_1,\dots,a_l)$
and $\mathbf{b}=(b_1,\dots,b_l)$ we put
\begin{eqnarray*}
\eta(\mathbf{a},\mathbf{b}) &=& \zeta(\underbrace{1,1,\dots,1}_{a_1-1},b_1+1,
\underbrace{1,1,\dots,1}_{a_2-1},b_2+1, \dots,
\underbrace{1,1,\dots,1}_{a_l-1},b_l+1) \\
|\mathbf{a}| &=& a_1 + a_2 + \dots + a_l \\
\binom{\mathbf{a} }{\mathbf{b} } &=& \binom{a_1}{b_1}\binom{a_2}{b_2}
\dots   \binom{a_l}{b_l} \\
(A,B)^{(\mathbf{a},\mathbf{b})} &=& A^{a_1} B^{b_1} \dots A^{a_l} B^{b_l}.
\end{eqnarray*}
Equation \eqref{eq.PhiKZ} implies that the KZ associator $\Phi^{\rm{KZ}}$
is Gevrey-0.
\end{proof}

\begin{remark}
\lbl{rem.associator}
It is not true that every associator $\Phi \in \A(\uparrow_3)$ is Gevrey-0.
In fact, it is not even true that the twist of a Gevrey-0 associator is
Gevrey-0, since the twist may have arbitrarily large coefficients.
\end{remark}

\begin{remark}
\lbl{rem.anotherproof}
There is an alternative proof of Theorem \ref{206} that does not use
associators. First decompose $T$ into smaller tangles, where
each smaller one is either elementary of type 1, or a braid. By
deformation we can assume that in any braid $X$, the horizontal
distance between any 2 strands is bigger than 1. Then the very
Kontsevich integral formula of $Z(X)$, see \cite{Ko} and
\cite[Sec.4.3]{B-N1} (and also \cite[Fig.13]{B-N1}, is regular and
easily seen to be Gevrey-0.
\end{remark}

\subsection{The LMO invariant}
\lbl{sub.LMOinv}

In this section we review the Aarhus version of the LMO invariant from
\cite[Part II]{BGRT}. For an equality of the Aarhus integral with the LMO
invariant, see \cite[Part III]{BGRT}.

We define a bilinear map

\begin{equation}
\lbl{eq.bracket}
\la \cdot, \cdot \ra: \A(\star_m) \otimes \A^{\tl}(\star_m) \to \A(\emptyset)
\end{equation}
as follows.  Suppose $\Gamma_1
\in \A(\star_m) $ and $\Gamma_2\in \A^{\tl}(\star_m)$ are Jacobi diagrams with
respectively  $k_{j}, l_j$ legs of label $j$, $j=1,2,\dots, m$. If
there is a $j$ such that $k_j \neq l_j$, let $\la \Gamma_1,
\Gamma_2\ra =0$, otherwise let $\la \Gamma_1, \Gamma_2\ra $ be the
sum of all possible ways to glue legs of label $j$ in $\Gamma_1$ to
legs of the same labels in $\Gamma_2$. Note that there are
$\prod_{j=1}^m (k_j)!$ terms in the sum. If $v\in \A(\star_m)$ and $u\in
\A^{\tl}(\star_m)$ have $k_j$ legs of label $j$, then

\begin{equation}
\lbl{200}
| \la v, u\ra| \le  \left (\prod_{j=1}^m k_{j}!\right) |v|\, |u|.
\end{equation}

It is known that the \ihs\ $M$ can be obtained
from $S^3$ by surgery along a {\em boundary link} $L$, where the
framing $\ve_1,\dots,\ve_m$ of the link components are $\pm 1$.
Suppose furthermore $L$ is the closure of a framed boundary string
link $T$. It is known \cite{LM} that

$$
Z_L = p  (\left [Z_T\right] \,\# \,
\left[\chi^{-1} (\Delta^{(m)}(\nu))\right]).
$$
Let us introduce some convenient notation. For Jacobi diagrams
$\Gamma_j\in \A(\star_1)$, $j=1,2,\dots,m$ let $\Gamma_1\otimes
\dots \otimes \Gamma_m \in \A(\star_m)$ be the union of all
$\Gamma_j$, with the legs of $\Gamma_j$  relabeled by $j$. Using
linearity we can define $v_1\otimes \dots \otimes v_m \in
\A(\star_m)$ for $v_j \in \A(\star_1)$. With the above notation, let
us define

\begin{eqnarray*}
\check Z_T &:=&  [Z_T] \,\#\,
\left[ \chi^{-1} (\Delta^{(m)}(\nu))\right]
\, \#\, \left[ \chi^{-1}(\nu^{\otimes m})\right], \\
\tilde Z_T &:=& \chi(\check Z_T) \sqcup E, \\
E &=& E(\ve_1,\dots,\ve_m) :=
\exp\left(-\frac{1}{2} \sum_{j=1}^m \e_j \strutn{j}{j} \right)
\in \A(\star_m).
\end{eqnarray*}
Notice that $\tilde Z_T$ has no struts.
Since $T$ is a boundary framed link, it follows from
\cite{HM} (see also \cite{GL0}) that $\tilde Z_T$ is treeless.
One can define

$$
\int T\, :=  \la   E, \tilde Z_T \ra  \in \A(\emptyset).
$$

Note that our  $\int T$ is  equal to $\int^{FG} \check Z_L$ in
\cite{BGRT}.

Suppose $U_\pm$ are the trivial string knot with framing $\pm 1$.
Suppose among $\ve_1,\ve_2,\dots,\ve_m$ there are $s_+$ positive
numbers and $s_-$ negative numbers. Then the LMO invariant of $M$
can be calculated by

$$
Z_M = \frac{\int T}{ (\int U_+)^{s_+} \, (\int U_-)^{s_-}}.
$$

\subsection{Proof of Theorem \ref{thm.gLMO}}
\lbl{sub.proofgLMO}

Using the multiplicative property of Gevrey-1 series, see Lemma
\ref{205}, to prove that $Z_M$ is Gevrey-1 it is enough to prove
the following lemma.

\begin{lemma}
\lbl{lem.boundarylemma}
For every boundary string link $T$ with framing
$\ve_1,\dots,\ve_m \in \{\pm 1\}$, the series $\int T$ is Gevrey-1.
\end{lemma}

\begin{proof}
By definition, we have
$$
E =\sum_{k_1,\dots,k_m \ge 0} E_{k_1,\dots,k_m}, \quad
$$
where

\begin{equation}
\lbl{209}
E_{k_1,\dots,k_m} := \prod_{j=1}^m \left ( \frac{(-\ve_j
/2)^{k_j}}{k_j!}\right) \left(\strutn{1}{1}\right)^{k_1}
\, \left(\strutn{2}{2}\right)^{k_2}
\dots \left(\strutn{m}{m}\right)^{k_m}
\end{equation}

Let $\tZ_T^{(2n,2k)}$ be the part of $\tZ_T$ of i-degree $2n$ and
e-degree $2k$. Since $\tZ_T$ is treeless, by Lemma \ref{101}, we
have $k \le n$, and hence the degree of $\tZ_T^{(2n,2k)}$ is less
than or equal to $2n$. By Theorem \ref{206}, Corollary \ref{207},
Proposition  \ref{204}(b), and the multiplicative property of Gevrey-0
series, $\tZ_T$ is Gevrey-0.  Thus, there is a constant $C$ such that for
every $n \ge 1$ we have:

$$
|\tZ_T^{(2n,2k)}| < C^n.
$$

Since $E$ consists of struts only,  $\la E,v\ra$ has degree equal
half the i-degree of $v$. Hence

\begin{equation}
\Grad_n \int T = \sum_{k=0}^n \la E, \tZ_T^{(2n,2k)}\ra \lbl{202}
\end{equation}

Recall that $E=\sum E_{k_1,\dots,k_m}$, and $E_{k_1,\dots k_m}$ has
$2k_j$ legs of label $j$.  For fixed $k$, the inner product $\la
E_{k_1,\dots,k_m}, \tZ_T^{(2n,2k)}\ra$ is non-zero only when
$k_1+\dots+k_m=k$. Using \eqref{200} and \eqref{209} we have

$$
\begin{aligned} |\la E_{k_1,\dots,k_m}, \tZ_T^{(2n,2k)}\ra | &< C^n \,
\prod_{j=1}^m \frac{(2k_j)!}{k_j!}\\
& < C^n \prod_{j=1}^m  2^{k_j} {k_j!} \qquad \text {since} \quad
\frac{(2k)!}{k!} \le 2^k \, k!\\
&< C^n 2^n n! \qquad \text {since}\quad  k_1+\dots k_m =k \le n
\end{aligned}
$$
It follows that the norm of

$$
\la E, \tZ_T^{(2n,2k)}\ra = \sum_{k_1 + \dots +k_m =k}
\la E_{k_1,\dots,k_m}, \tZ_T^{(2n,2k)}\ra
$$
can be estimated by

$$
\la E, \tZ_T^{(2n,2k)}\ra < C^n 2^n n! \left(\,  \sum_{k_1 + \dots k_m =k}
1 \right)= C^n 2^n n!\,  \binom{k+m-1}{m} \le C^n 2^n n! 2^{n+m}.
$$
Using \eqref{202}, we get

$$
| \Grad_n \int T| < n C^n 2^n n! 2^{n+m} < n! \, C'^n
$$
for an appropriate constant $C'=C'_T$. This concludes the proof of Theorem
\ref{thm.gLMO}.
\end{proof}

\begin{remark}
\lbl{rem.qhs}
Without doubt, Theorem \ref{thm.gLMO} holds for rational homology spheres
as well. This requires a technical modification of the proof that allows
one to deal with Jacobi diagrams with tree components. This is possible,
but it requires another layer of technicalities that we will not present
here.
\end{remark}

\subsection{Proof of Theorem \ref{thm.4}}
\lbl{sub.proofthm4}

Theorem \ref{thm.4} follows immediately from Theorem \ref{thm.gLMO} and
the following Lemma \ref{lem.weightbound}.

\begin{lemma}
\lbl{lem.weightbound} For every simple Lie algebra $\fg$ there is a
constant $C$ such that for any  Jacobi diagram $\Ga \in
\A(\emptyset)$ of degree $n>0$ we have:
$$
|W_{\fg}(\Ga)| \leq  C^n
$$
\end{lemma}


\begin{proof}
 $\Ga$
is obtained from a cloud of $2n$ $Y$ graphs by a complete pairing of their
legs. By the definition of the weight system, it follows that $W_{\fg}(\Ga)$
is obtained by the contraction of the indices of a tensor in $\oplus_{2n}
(\fg^{\otimes 3})$. The result follows.
\end{proof}

\subsection{On the dimension of the space of Feynman diagrams}
\lbl{sub.dimension}

Although our proof of Theorem \ref{thm.gLMO} is completed, in this section
we will give some estimates for the space $\calD_n(\emptyset)$ and
$\calA_n(\emptyset)$ of Feynman diagrams, introduced in Section \ref{sub.LMO}.
Our crude estimates explain the Gevrey nature of the Gromov norm of the LMO
invariant.

Let $\calD_n(\emptyset)$ denote the (finite dimensional) vector space
with basis the set of Jacobi diagrams with no legs of degree $n$,
and $\calA_n(\emptyset)$ is its quotient by the AS and IHX relations.

Let us say that a  Jacobi diagram is {\em normalized} if it is made out of a
number of
disjoint circles together with  a number of chords
that each begin and end on the same circle. Let $\calS$ (resp. $\calS_n$)
denote the set of normalized graphs (resp. of degree $n$).

\begin{lemma}
\lbl{lem.normalized}
\rm{(a)} $\dim \calD_n(\emptyset) \leq n!^3 C'^n$ for some $C'$.
\newline
\rm{(b)} $\calS_n$ is a  spanning set for $\calA_n(\emptyset)$.
\newline
\rm{(c)} $ \dim \calA_n(\emptyset) \leq n! C^n$ for some $C$.
\end{lemma}

\begin{proof}
If $\Ga$ is a trivalent graph of degree $n$, then we can cut it along each
of its edges. We obtain a cloud of $2n$ $Y$ graphs. $\Ga$ can be reconstructed
by matching the legs of the $Y$ graphs. There are $6n$ legs, and they can
be matched in $(6n)!!$ ways. Using
Stirling's formula \cite[p.50-53]{F}

\begin{equation}
\lbl{eq.stirlingformula}
\sqrt{2\pi} \,\, n^{n+1/2} e^{-n + 1/(12n+1)} < n!
\sqrt{2\pi} \,\, n^{n+1/2} e^{-n + 1/(12n)}.
\end{equation}
the result follows. This proves (a). The above bound may feel a little
crude, since we did not take into account automorphisms of the graphs.
Nevertheless, it seems to be asymptotically optimal; see also
\cite[p.55.Cor.2.17]{Bo}.

For parts (b,c), we need to understand what we gain by the AS and IHX relation.
If $\Ga$ is a connected trivalent graph of degree $n$, choose a cycle in
it. Then, using the IHX relation repeatedly, write $\Ga$ as a
linear combination of ``chord diagrams'' on that cycle.
Since there are at most $(2n)!!$ chord diagrams with $n$ chords on a circle,
the result follows for connected graphs. Applying the above reasoning
to each connected component of a trivalent graph implies the result in
general.
\end{proof}

\subsection{An integral version of the  Ohtsuki series}
\lbl{sub.ohtsukiseries}

In quantum topology, there are two commonly used Taylor series
expansions of an element $f(q)$ of the Habiro ring; namely setting
$q=e^{1/x}$ or setting $q=1+1/x$. So far we have worked with $q=
e^{1/x}$. The other substitution $q= 1+ 1/x$ leads to another map

\begin{equation*}
T^{\BZ}: \Lhat \longto \BZ[[1/x]], \qquad (T^{\BZ}f)(x)=f(1+1/x)
\end{equation*}
which is also injective. We may also consider a map:
\begin{equation*}
F^{\BZ}:  \text{Knotted Objects} \longto \BZ[[1/x]],
\qquad
F^{\BZ}=T^{\BZ} \circ \Phi.
\end{equation*}

The interest in the latter formal power series lies in the fact that
it has integer coefficients. In fact, the original definition of the
Ohtsuki series is in this form, see \cite{Oh1}.

From the point of view of analysis, the series $F_{\calK}(x)$ and
$F^{\BZ}_{\calK}(x)$
are simple reparametrizations of one another, by an analytic change of
variables.
Our next lemma shows that the notion of a Gevrey series is independent
of an analytic change of variables.

\begin{lemma}
\lbl{lem.equiv}
Consider a formal power series $f(x) \in \BC[[1/x]]$ and let
$g(x)=f(e^{1/x}-1) \in \BC[[x]]$. Then $f(x)$ is Gevrey iff $g(x)$ is Gevrey.
\end{lemma}

\begin{proof}
Let
\begin{eqnarray}
\lbl{eq.f1}
f(x) &= & \sum_{k=0}^\infty a_k \frac{1}{x^k} \\
\lbl{eq.f2}
f(e^{1/x}-1) &= & \sum_{k=0}^\infty b_k \frac{1}{x^k}
\end{eqnarray}
Then the sequences $(a_n)$ and $(b_n)$ are related by an upper-triangular
matrix with $1$ on the diagonal.
The asymptotic behavior of the entries of this matrix make the lemma
possible. For a thorough discussion on that subject, see also
Hardy's book \cite{Ha}.
The entries of the matrix are given by Stirling numbers.
The {\em Stirling numbers} $s_{n,k}$ of the first kind satisfy:
$$
\frac{1}{x^k}=k! \sum_{n=k}^\infty  \frac{s_{n,k}}{n!} (1-e^{-1/x})^n.
$$
Substituting for $1/x^k$ from the above identity into \eqref{eq.f2}
and rearranging, it follows that:


\begin{eqnarray*}
a_n 
& = & \sum_{k=0}^n b_{n-k} (-1)^{k} \frac{(n-k)!}{n!} s_{n,n-k}.
\end{eqnarray*}
Suppose now that $(b_n)$ is Gevrey:
$$
|b_n| \leq n! C' C^n.
$$
On the other hand, we have:
$$
s_{n,n-k} =\frac{n^{2k}}{2^k k!}\left(1+\frac{c_1(k)}{n} + \frac{c_2(k)}{n^2}
+ \dots \right)
$$
where $c_1(k)=-k(2k+1)/3, \dots$ which are polynomials
in $k$.

Since $s_{n,n-k} \geq 0$ for all $n$ and $k$ with
$ k \leq n$, and
$$
\left(\frac{(n-k)!}{n!}\right)^2 n^{2k}=
\left( \frac{n^k}{(n-k+1)\dots n} \right)^2 \leq 1
$$
it follows that
\begin{eqnarray*}
|a_n| & = & | \sum_{k=0}^n b_{n-k} (-1)^{k} \frac{(n-k)!}{n!} s_{n,n-k}| \\
& \leq & n! C' C^n \sum_{k=0}^n   C^{-k} \left(\frac{(n-k)!}{n!}\right)^2
\frac{n^{2k}}{2^k k!}
\left(1 + \frac{c_1(k)}{n} + \frac{c_2(k)}{n^2}
+ \dots \right) \\
&=& n! C' C^n \sum_{k=0}^\infty (2C)^{-k}\frac{1}{k!}
\left( 1 + \frac{c_1(k)}{n} + \frac{c_2(k)}{n^2} + \dots
\right) \\
&=& n! C' C^n e^{-(2C)^{-1}} \left(1 +
\frac{d_1(k)}{n} + \frac{d_2(k)}{n^2} + \dots
\right),
\end{eqnarray*}
where $d_1(k), d_2(k), \dots$ are polynomials in $k$.
This concludes one half of the theorem. The other half is similar.
\end{proof}

\begin{remark}
\lbl{rem.cg}
In \cite{CG3} a more general statement is shown. Namely, suppose that
$f(x) \in \BC[[1/x]]$ is a power series and $\tau(x)$ is analytic in a
neighborhood of infinity and small (i.e., $\tau(\infty)=0$). Consider the
power series $f^{\tau}(x)=f(1/x+\tau(x))$. Then, $f(x)$ is Gevrey iff
$f^{\tau}(x)$ is Gevrey.
\end{remark}

Theorem \ref{thm.2} and Lemma \ref{lem.equiv} imply that:

\begin{corollary}
\lbl{cor.ohtsuki} For every knotted object $\calK$, the integral
Ohtsuki series $F^{\BZ}_{\calK}(x) \in \BZ[[1/x]]$ is Gevrey.
\end{corollary}

\section{Proof of Theorem \ref{thm.2}}
\lbl{sec.main}

In this section we give a proof of Theorem \ref{thm.2}.
To simplify notation, let $\la f(x) \ra_k$ denote the coefficient of $1/x^k$
in a formal power series $f(x)$.

\begin{lemma}
\lbl{lem.estimate}
\rm{(a)}
If two sequences $(f_n(q))$ and $(g_n(q))$ are nicely bounded, so is their
product.
\newline
\rm{(b)}
If a sequence $(f_n(q))$ is nicely bounded
then there exist constants $C,C'$ such that
for every $k$ and every $n$ we have:
$$
| \la f_k(e^{1/x}) \ra_n | \leq \frac{1}{n!} k^{2 n} e^{C' n + C k}.
$$
In particular, there exist constants $C''$ such that
for every $n$ and every $k$ with $0 \leq k \leq n$ we have:
$$
| \la f_k(e^{1/x}) \ra_n | \leq n! e^{C'' n}
$$
\rm{(c)} The sequence $((q)_n)$ is nicely bounded.
\end{lemma}

\begin{proof}
Part (a) is easy.

For part (b), we may write
$$
f_k(q)=\sum_{j } a_{k,j} q^j.
$$
$j \in [C' k^2 + c', C'' k^2 + c'']$, and
$|a_{k,j}| \leq e^{C k}$ for all such $j$.  It follows that
\begin{eqnarray*}
| \la  f_k(e^{1/x}) \ra_n| &=&
| \sum_{j } a_{k,j} \la e^{j/x} \ra_n | \\
& = & \frac{1}{n!}
| \sum_{j } a_{k,j} j^{n}| \\
& \leq & \frac{1}{n!} e^{C k}
\sum_{j } j^{n} \\
&=& \frac{1}{n!} e^{C k}
k^{2n} e^{C' n}.
\end{eqnarray*}
If in addition $k \leq n$, then Stirling's formula
\eqref{eq.stirlingformula} and the above implies that
\begin{eqnarray*}
| \la  f_k(e^{1/x}) \ra_n| & \leq & \frac{1}{n!} e^{C k}
k^{2n} e^{C' n} \\
& \leq & \frac{1}{n!} e^{C n}
n^{2n} e^{C' n} \\
& \leq & n! e^{C'' n}.
\end{eqnarray*}

For part (c), it is easy to see that
\begin{eqnarray*}
\span_q (q)_n &=& [0,n(n+1)/2] \\
||(q)_n||_1 & \leq & 2^n.
\end{eqnarray*}
\end{proof}

\begin{proof}(of Theorem \ref{thm.2})
Let us fix an element
$$
f(q)=\sum_{n=0}^\infty f_n (q) (q)_n
$$
of $\Lhat^b$ where $(f_n(q))$ is nicely bounded. Since
$\la (e^{1/x})_k \ra_n=0$ for $k >n$, it follows that
\begin{equation}
\lbl{eq.fff}
\la f(e^{1/x}) \ra_n=\sum_{k=0}^n \la f_k(e^{1/x}) (e^{1/x})_k
\ra_n.
\end{equation}
Lemma \ref{lem.estimate} (a) and (c) implies that the sequence $(f_n(q)(q)_n)$
is nicely bounded, and therefore by (b) there exist a constant so that
for every $n$ and every $k$ with $0 \leq k \leq n$ we have:
\begin{eqnarray*}
| \la f_k(e^{1/x})(e^{1/x};e^{1/x})_k) \ra_n |
& \leq & n! C''^n.
\end{eqnarray*}
Using Equation \eqref{eq.fff}, the result follows.
\end{proof}

\section{Proof of Theorem \ref{thm.1}}
\lbl{sec.proofthm1}

In this section we give a proof of Theorem \ref{thm.1} using two
properties of the colored Jones polynomial: $q$-holonomicity (see
\cite{GL1}), and integrality, due to Habiro \cite{H1}.

We call a {\em virtual $sl_2$-module} any $\BQ(q^{1/4})$-linear
combination of finite-dimensional $sl_2$-modules. Suppose $L$ is
framed oriented link with $m$ numbered components, and
$U_1,\dots,U_m$ are virtual $sl_2$-modules, then there is defined
the colored Jones polynomial (rather rational function)
$$
J_L(U_1,\dots,U_r) \in \BQ(q^{\pm/4})
$$
normalized by the quantum dimension for the unknot; see \cite{Tu}.
Habiro introduced two important sequences of virtual $sl_2$, denoted
by $P'_n$ and $P''_n$, see \cite{H1}.

\subsection{Proof of Theorem \ref{thm.1} for knots}
\lbl{sub.newb.knots}

Let us fix a knot $K$ with framing 0 in 3-space. \cite{HL} give the
following formula for the Kashaev invariant $\Phi_K(q)$:

\begin{equation}
\lbl{eq.kknot}
\Phi_K(q)=\sum_{n=0}^\infty J_K(P''_n)(q) (q)_n(q^{-1})_n
\in \Lhat
\end{equation}
where, according to Habiro we have: $J_K(P''_n)(q) \in \BZ[q^{\pm 1}]$
for all $n$; see \cite{H2}.

In \cite{GL1} we proved that the sequence $(J_K(P''_n))$ is $q$-holonomic
and in \cite{GL2} we proved that it is also nicely bounded. It follows
that the sequence $ (J_K(P''_n)(q^{-1})_n)$ is $q$-holonomic and
nicely bounded (by Lemma \ref{lem.estimate}). Thus, $\Phi_K(q)$
lies in $\Lhathol \cap \Lhat^b$. This concludes the proof of Theorem
\ref{thm.1} for knots.

Although we will not need it here, let us mention that \cite{HL}
prove that when $q=e^{2 \pi i/N}$, then
$$
\Phi_K(e^{2 \pi i/N})=\la K \ra_N
$$
where $\la K \ra_N$ is the well-known {\em Kashaev invariant} of
a knot $K$; see \cite{Ka} and \cite{MM}.

\subsection{Proof of Theorem \ref{thm.1} for \ihs s}
\lbl{sub.newb.ihs}

\begin{proof}(of Theorem \ref{thm.1})
The proof will use integrality and holonomicity properties
of the colored Jones function of a link in 3-space.

Consider an \ihs\ $M$. We can find surgery presentation
$M=S^3_{L,f}$ where $L$ is an algebraically split
link $L$ in $S^3$ of $r$ ordered components, with framing
$f=(f_1,\dots, f_r)$.

Habiro considers the following series:
\begin{equation}
\lbl{eq.hab1}
\Phi_M(q)= \sum_{k_1,\dots, k_r=0}^\infty
J_L(P'_{k_1},\dots, P'_{k_r}) \prod_{i=1}^r (-f_i q)^{-f_i k_i(k_3+3)/4}.
\end{equation}
where $J_L$ is the colored Jones function of the $0$-framed link $L$.
Habiro proves that
\begin{itemize}
\item
for all $k_1,\dots, k_r \in \BN$, we have:
\begin{equation}
\lbl{eq.hab2}
J_L(P'_{k_1},\dots, P'_{k_r}) \in
\frac{\{2m+1\}!}{\{m\}!\{1\}} \,\,\BZ[q^{\pm1/2}],
\end{equation}
where $m=\max\{k_1,\dots, k_r\}$, and
\begin{equation}
\lbl{eq.aa}
\{a\}!:=\prod_{j=1}^a(q^{a/2}-q^{-a/2})=(-1)^a q^{-a(a+1)/2} (q)_a
\end{equation}
Moreover,
\begin{equation}
\lbl{eq.integ}
J_L(P'_{k_1},\dots, P'_{k_r}) \prod_{i=1}^r (-f_i q)^{-f_i k_i(k_3+3)/4}
\in \BZ[q^{\pm 1}]
\end{equation}
for all $k_1,\dots, k_r \in \BN$.
Thus, $\Phi_M(q) \in \Lhat$ is a convergent series.
\item
The right hand side of Equation \eqref{eq.hab1} is independent of the surgery
presentation $M=S^3_{L,f}$, and depends on $M$ alone.
\end{itemize}

To simplify notation, let us define:
\begin{eqnarray}
\lbl{eq.aL}
\a(k_1,\dots,k_r) &=& \frac{1}{(q)_m} J_L(P'_{k_1},\dots, P'_{k_r})
\prod_{i=1}^r (-f_i q)^{-f_i k_i(k_3+3)/4}  \\
\lbl{eq.bL}
\b(m) &=& \sum_{k_1,\dots, k_r; \, \max\{k_1,\dots,k_r\}=m}
\a(k_1,\dots,k_r)
\end{eqnarray}

Since $\frac{\{2m+1\}!}{\{m\}!\{1\}}$ is divisible by $\{m\}!$, it follows that
for all $k_1,\dots,k_r$, we have
$\a(k_1,\dots,k_r) \in \BZ[q^{\pm 1}]$, and consequently, for all $m$
we have $\b(m) \in \BZ[q^{\pm 1}]$. Moreover,

\begin{equation}
\lbl{eq.newIM}
\Phi_M(q)=\sum_{m=0}^\infty \b(m) (q)_m
\end{equation}

Theorem \ref{thm.1} for \ihs s
follows from the following

\begin{theorem}
\lbl{thm.newbb}
For every unit-framed algebraically split link $(L,f)$, the sequence
$(\a(m))$ is $q$-holonomic and nicely bounded.
\end{theorem}

\begin{proof}
It suffices to show that $(\b(m))$ is $q$-holonomic, and nicely
bounded.

Let us first recall that the class of $q$-holonomic functions in several
variables  is closed under the operations of
\begin{itemize}
\item[(P1)] sum,
\item[(P2)] product,
\item[(P3)] specialization,
\item[(P4)] definite summation
\item[(P5)] contains the proper $q$-hypergeometric functions.
\end{itemize}
For a proof, see  \cite{Ze}.

Without loss of generality, let us assume that $r=2$ (the general case
follows from inclusion-exclusion). Then, we have:

\begin{equation}
\lbl{eq.bbm}
\b(m)=\sum_{k_1=0}^m \a(k_1,m) + \sum_{k_2=0}^m \a(m,k_2) - \a(m,m).
\end{equation}

Changing basis from $\{P'_k\}$ to $\{V_l \}$ it follows that
$$
\{k_1\}! \{k_2 \}!   J_L(P'_{k_1}, P'_{k_2})=\sum_{l_1=0}^{k_1}
\sum_{l_2=0}^{k_2} P_{k_1,k_2}^{l_1,l_2} J_L(V_{l_1},V_{l_2}),
$$
where $P_{k_1,k_2}^{l_1,l_2} \in \BZ[q^{\pm/2}]$ are explicit Laurent
polynomials which are proper $q$-hypergeometric; see \cite[Sec.4]{GL2}.
This, together with Equation \eqref{eq.aL} implies that
\begin{equation}
\lbl{eq.jj}
\{m \}! \{k_1\}! \{k_2 \}! \a(k_1,k_2)=
\sum_{l_1=0}^{k_1}
\sum_{l_2=0}^{k_2} R_{k_1,k_2}^{l_1,l_2} J_L(V_{l_1},V_{l_2})
\end{equation}
where $R_{k_1,k_2}^{l_1,l_2}(q) \in \BZ[q^{\pm/2}]$ are
proper $q$-hypergeometric Laurent polynomials, and $m=\max\{k_1,k_2\}$.

Now $J_L(V_{l_1},V_{l_2})$ can be written as a multisum:

\begin{equation}
\lbl{eq.cj}
J_L(V_{l_1},V_{l_2})=\sum_{j_1=0}^{l_1} \sum_{j_2=0}^{l_2}
F_{l_1,l_2,j_1,j_2}
\end{equation}
where $F_{l_1,l_2,j_1,j_2} \in \BZ[q^{\pm/2}]$ is a proper
$q$-hypergeometric summand; see \cite[Sec.3]{GL1}.

Equations \eqref{eq.jj}, \eqref{eq.cj} and Properties P4, P5 imply that
$\a(k_1,k_2)$ is $q$-holonomic in both variables $(k_1,k_2)$.
Together with Property P4, it follows that
$\sum_{k_1=0}^r \a(k_1,s)$ is $q$-holonomic in both variables
$(r,s)$, and (by Property P3)
$\sum_{k_1=0}^m \a(k_1,m)$ is $q$-holonomic in $m$. Alternatively, the
WZ algorithm of \cite{WZ}
and Equations \eqref{eq.jj} and \eqref{eq.cj} imply directly
(and constructively) that $\sum_{k_1=0}^m \a(k_1,m)$ is $q$-holonomic in $m$.

Likewise, $\sum_{k_2=0}^m \a(m,k_2)$ is $q$-holonomic in $m$, and $\a(m,m)$ is
$q$-holonomic in $m$. Property P1 and Equation \eqref{eq.bbm}
implies that $\b(m)$ is $q$-holonomic in $m$.

Alternatively, we could have used the identity

\begin{equation}
\lbl{eq.balt}
\b(m)=\sum_{0 \leq k_1,\dots, k_r \leq m}\a(k_1,\dots,k_r)-
\sum_{0 \leq k_1,\dots, k_r \leq m-1}\a(k_1,\dots,k_r)
\end{equation}
and the $q$-holonomicity of $\a(k_1,\dots,k_r)$ (as follows by the WZ
algorithm) to deduce the $q$-holonomicity of $\b(m)$.

It remains to show that $\b(m)$ is nicely bounded.
Let us say that a multi-indexed sequence $(f(r_1,r_2,\dots))$
of Laurent polynomials is {\em nicely bounded} if it satisfies
\eqref{eq.degbound} and \eqref{eq.coeffbound} for all $r_1,r_2,\dots$
with $r_1, r_2, \dots \leq n$.

It is easy to see that the class of nicely bounded functions satisfies
properties P1-P4 and contains the proper $q$-hypergeometric terms that are
Laurent polynomials.

Repeating our previous steps, Equations \eqref{eq.cj} and \eqref{eq.jj}
imply that $\{m \}! \{k_1\}! \{k_2 \}! \a(k_1,k_2)$
is nicely bounded as a function of both variables $(k_1,k_2)$.
Lemma \ref{lem.boyd} below (communicated to us by D. Boyd), implies that
$\a(k_1,k_2)$ is nicely bounded as a function of both variables
$(k_1,k_2)$.

Our previous steps (or Equation \eqref{eq.balt})
now imply that $\b(m)$ is nicely bounded.
\end{proof}

This concludes the proof of Theorem \ref{thm.1} for integer homology
spheres.
\end{proof}

\begin{lemma}
\lbl{lem.boyd}(Boyd)
If $(f_n(q))$ is a sequence of Laurent polynomials such that
$((q)_n f_n(q))$ is nicely bounded, then $(f_n(q))$ is nicely bounded, too.
\end{lemma}

For a proof, see \cite[Sec.7]{GL2}.

\begin{remark}
\lbl{rem.special}
In the special case where $M$ is obtained by $\pm 1$ surgery on a knot in
3-space, Lawerence-Ron have shown independently
that the formal power series $F_M(x)$ is Gevrey;
see \cite{LR}.
\end{remark}

\ifx\undefined\bysame
        \newcommand{\bysame}{\leavevmode\hbox
to3em{\hrulefill}\,}
\fi


\begin{thebibliography}{[EMSS]}

\bibitem[An]{An} G.E. Andrews,
        {\em The theory of partitions},
        Cambridge University Press, Cambridge, 1998.

\bibitem[Ba]{Ba} W. Balser,
        {\em From divergent power series to analytic functions.
        Theory and application of multisummable power series},
        Lecture Notes in Mathematics, {\bf 1582} Springer-Verlag, Berlin,
        1994.

\bibitem[Bo]{Bo} B. Bollob\'as,
        {\em Random graphs},
        Second edition, Cambridge Studies in Advanced Mathematics, {\bf 73}
        Cambridge University Press 2001.

\bibitem[B-N1]{B-N1} D. Bar-Natan,
        {\em On the Vassiliev knot invariants},
        Topology {\bf 34} (1995)  423--472.

\bibitem[B-N2]{B-N2} \bysame,
        {\em Non-Associative Tangles},
        In: Geometric Topology (W. Kazez, Ed.), Proc. Georgia Int.
        Topology  Conf. 1993
        AMS/IP Studies in Advanced Mathematics {\bf 1997} 139--183.

\bibitem[B-N3]{B-N3} \bysame,
        {\em On associators and the Grothendieck-Teichmuller group I},
        Selecta Math.   {\bf 4}  (1998) 183--212.

\bibitem[BGRT]{BGRT} \bysame, S. Garoufalidis, L. Rozansky and
        D. Thurston,
        {\em The Aarhus integral of rational homology 3-spheres I-III},
        Selecta Math. {\bf 8} (2002) 315--339,
        {\bf 8} (2002) 341--371, {\bf 10} (2004) 305--324.

\bibitem[BLT]{BLT}
        \bysame, T.~T.~Q.~L\^e, and D.~P.~Thurston,
        {\em Two applications of elementary knot theory to Lie algebras
        and Vassiliev invariants},
        Geometry and Topology {\bf 7} (2003), 1-31

\bibitem[CG1]{CG1} O. Costin  and S. Garoufalidis,
        {\em Resurgence of the Kontsevich-Zagier power series},
        preprint 2006 {\tt math.GT/0609619}.

\bibitem[CG2]{CG2} \bysame and \bysame,
        {\em Resurgence of the Euler-MacLaurin summation formula},
        preprint 2006.

\bibitem[CG3]{CG3} \bysame and \bysame,
        {\em Resurgence of 1-dimensional sums of $q$-factorials},
        preprint 2007.

\bibitem[Dr]{Dr} V. G. Drinfeld,
        {\em On quasitriangular quasi Hopf algebras and a group
        closely connected with $Gal(\overline{\mathbb Q}/\mathbb Q)$},
        Leningrad Math. J. {\bf 2} (1991) 829--860.

\bibitem[Ec]{Ec} J. \'Ecalle,
        {\em Resurgent functions},
        Vol. I-II, Mathematical Publications of Orsay {\bf 81}, 1981.

\bibitem[F]{F} W. Feller,
        {\em An Introduction to Probability Theory and Its Applications},
        Vol. 1, 3rd ed. New York, Wiley, 1968.

\bibitem[GL0]{GL0} S. Garoufalidis and J. Levine,
        {\em Tree-level invariants of 3-manifolds, Massey products and the
        Johnson homomorphism},
        Graphs and patterns in mathematics and theoretical physics,
         Proc. Sympos. Pure Math. AMS {\bf 73} (2005) 173--203.

\bibitem[GL1]{GL1} \bysame and  T.T.Q. Le,
        {\em The colored Jones function is $q$-holonomic}
        Geom. and Topology {\bf 9} (2005) 1253--1293.

\bibitem[GL2]{GL2} \bysame and \bysame,
        {\em Asymptotics of the colored Jones function},
        preprint 2005 {\tt math.GT/0508100}.

\bibitem[G]{G} \bysame,
        {\em Difference and differential equations for the colored Jones
        function},
        Journal of Knot Theory and its Ramifications, {\em in press}.

\bibitem[GG]{GG} \bysame and J. Geronimo,
        {\em  Asymptotics of $q$-difference equations},
        JAMI Conference Proceedings, in press.

\bibitem[HM]{HM} N. Habegger and G. masbaum,
        {\em The Kontsevich integral and Milnor's invariants},
        Topology  {\bf 39}  (2000) 1253--1289.

\bibitem[Ha]{Ha} G.H. Hardy,
        {\em Divergent Series},
        Oxford  Clarendon Press (1949).

\bibitem[H1]{H1} K. Habiro,
        {\em On the quantum $\mathfrak{sl}_2$
        invariants of knots and integral homology
        spheres},
        Geom. Topol. Monogr. {\bf 4} (2002) 55--68.

\bibitem[H2]{H2} \bysame,
        {\em Cyclotomic completions of polynomial rings},
        Publ. Res. Inst. Math. Sci.  {\bf 40}  (2004) 1127--1146.

\bibitem[H3]{H3} \bysame,
        {\em Bottom tangles and universal invariants},
        preprint 2005 {\tt math.GT/0505219}.

\bibitem[HL]{HL} V. Huynh and T.T.Q. Le,
        {\em  On the Colored Jones Polynomial and the Kashaev invariant},
        preprint 2005 {\tt math.GT/0503296}.


\bibitem[Ka]{Ka} R. Kashaev,
        {\em The hyperbolic volume of knots from the quantum dilogarithm},
        Modern Phys. Lett. A {\bf 39} (1997) 269--275.

\bibitem[Ko]{Ko} M.~Kontsevich,
        {\em  Vassiliev's knot invariants},
        Advances in Soviet Mathematics, {\bf 16} (1993), 137-150.

\bibitem[LR]{LR} R. Lawrence, O. Ron,
        {\em On Habiro's cyclotomic expansions of the Ohtsuki invariant},
        Journal of Knot Theory and its Ramifications, {\bf 15} (2006) 661--680.

\bibitem[Le1]{Le1} T. T. Q. L\^e,
        {\em Integrality and symmetry of quantum link invariants},
        Duke Math. J.  {\bf 102}  (2000) 273--306.

\bibitem[Le2]{Le5} \bysame,
        {\em Quantum invariants of 3-manifolds: integrality,
        splitting, and perturbative expansion},
        Topology Appl.  {\bf 127} (2003)  125--152.

\bibitem[LM]{LM} \bysame and  J.~Murakami,
        {\em The universal Vassiliev-Kontsevich invariant for framed
        oriented links},
        Compositio Math., {\bf 102} (1996)   41--64.

\bibitem[LMO]{LMO}  \bysame, J. Murakami and  T. Ohtsuki,
        {\em A universal quantum invariant of 3-manifolds},
        Topology {\bf 37} (1998) 539--574.

\bibitem[MM]{MM} H. Murakami and J. Murakami,
        {\em The colored Jones polynomials and the simplicial volume of a
        knot},
        Acta Math.  {\bf 186}  (2001) 85--104.

\bibitem[Oh1]{Oh1} T. Ohtsuki,
        {\em A polynomial invariant of rational homology $3$-spheres},
        Invent. Math. {\bf 123}  (1996) 241--257.

\bibitem[Oh2]{Ohtsukibook} \bysame,
        {\em Quantum invariants. A study of knots, 3-manifolds,
        and their sets},
        Series on Knots and Everything, {\bf 29},
        World Scientific Publishing Co., Inc. 2002.

\bibitem[Ra]{Ra} J.P. Ramis,
        {\em S\'eries divergentes et th\'eories asymptotiques},
        Bull. Soc. Math. France  {\bf 121}  (1993).


\bibitem[Tu]{Tu} V. Turaev,
        {\em The Yang-Baxter equation and invariants of links},
        Inventiones Math. {\bf 92} (1988) 527--553.

\bibitem[Vo]{Vo} P. Vogel,
        {\em Invariants de type fini},
        Nouveaux invariants en g\'eom\'etrie et en topologie,
        Panor. Synth\`eses {\bf 11} Soc. Math. France (2001) 99--128.

\bibitem[WZ]{WZ} H. Wilf and D. Zeilberger,
        {\em An algorithmic proof theory for hypergeometric (ordinary and
        $q$) multisum/integral identities},
        Inventiones Math. {\bf 108} (1992)  575--633.

\bibitem[Za]{Za} D. Zagier,
        {\em Vassiliev invariants and a strange identity related to the
        Dedekind eta-function},
        Topology  {\bf 40}  (2001) 945--960.

\bibitem[Ze]{Ze} D. Zeilberger,
        {\em A holonomic systems approach to special functions identities},
        J. Comput. Appl. Math. {\bf 32} (1990) 321--368.

\end{thebibliography}
\end{document}